\documentclass[11pt]{article}

\usepackage{amssymb,epsfig}
\usepackage{multicol}
\usepackage{tabularx}
\usepackage{multirow}

\newcommand{\pph}{\textsc{pph} }

\newcommand{\chaikin}{\textsc{c} }
\newcommand{\pphapprox}{\textsc{ppha} }

\newcommand{\R}{\mathbb{R}}
\newcommand{\Z}{\mathbb{Z}}
\newcommand{\demi}{\frac{1}{2}}

   \title {\bf  On a $C^2$-nonlinear subdivision scheme avoiding Gibbs oscillations}

\author{S. Amat \thanks{
 Departamento de Matem\'atica Aplicada y Estad\'{\i}stica.
   Universidad  Polit\'ecnica de Cartagena (Spain).
Research supported in part by the Spanish grants MTM2007-62945.
e-mail:{\tt sergio.amat@upct.es}} \and
  K. Dadourian\thanks{ Ecole Centrale de Marseille.
Laboratoire d'Analyse Topologie et Probabilites.
 e-mail:{\tt
dadouria@cmi.univ-mrs.fr}} \and J. Liandrat\thanks{ Ecole Centrale
de Marseille. Laboratoire d'Analyse Topologie et Probabilites.
 e-mail:{\tt
jliandrat@ec-marseille.fr}}}
\begin{document}

\maketitle

\begin{abstract}
   This paper is devoted to the presentation and the study of a new  nonlinear
   subdivision scheme eliminating the Gibbs oscillations close to discontinuities.
Its convergence, stability and
   order of approximation are analyzed. It is proved that this scheme converges
   towards limit functions of H\"older regu\-la\-rity index  larger  than $1.192$.
   Numerical estimates provide an  H\"older regu\-la\-rity index of $2.438$.
 Up to our knowledge, this subdivision scheme is the first one
  that achieves simultaneously  the control of the Gibbs phenomenon and
 regularity index larger than $1$ for its limit functions.
    \end{abstract}

{\bf Key Words.}  Nonlinear subdivision scheme, limit function,
regularity, stability, Gibbs phenomenon.

\vspace{10pt} {\bf AMS(MOS) subject classifications.} 41A05,
41A10, 65D05, 65D17

\newtheorem{lemma}{Lemma}
\newtheorem{proof}{Proof}
\newtheorem{definition}{Definition}
\newtheorem{remark}{Remark}
\newtheorem{proposition}{Proposition}
\newtheorem{theorem}{Theorem}

    \section{Introduction}

Subdivision schemes are useful tools for generating smooth curves
and surfaces.  For convergent schemes, starting from discrete sets
of control points and using basic rules of low complexity,  curves
or surfaces  can be obtained as  limits (called limit functions)
of  sequences of points generated by recursive applications of the
subdivision scheme.

A simple example of subdivision scheme is the family of interpolatory subdivision schemes,
 based on Lagrange's interpolation that have been  derived and analyzed in \cite{DeDu}.
 Other example is the family of   spline subdivision schemes related to spline
 spaces \cite{CC78}.

The four-point interpolatory scheme \cite{Du}, \cite{DyGe} is a
convergent linear scheme of the first family, involving four-point
stencils  at each subdivision, for which the limit function is at
least in the space\footnote{For $0<\alpha<1$, $f \in C^\alpha(\R)$
iff $f$ is bounded and $\exists C>0$ such that $\forall x,y \in
\R, |f(x)-f(y)|\leq C|x-y|^\alpha$\\ For $\alpha>1$,  $f \in
C^\alpha(\R)$ iff $f^{(\left[\alpha\right])}$ is bounded and
$f^{(\left[\alpha\right])} \in C^{(\alpha-\left[\alpha\right])}$
where  $\left[\alpha\right]$ is the integer part of $\alpha$}
$C^1$. The Chaikin algorithm \cite{chaikin} is an example of
spline subdivision scheme, with lower complexity as the previous
example and converging towards  $C^{2-}$ functions\footnote{
$C^{\alpha-}=\{f \in C^\beta, \forall \beta <\alpha\}$ }.

For applications, for instance to computer aided geometric  design
or image processing, complexity and convergence/regularity are not
the only quality criterion. The order of approximation, that characterizes
the precision of the scheme, is an other important one. Moreover,
oscillations that could occur in the limit function at the
vicinity of strongly variating data (coming from the sampling of discontinuous functions) , called Gibbs oscillation, are
really undesirable.

In the last decade, various attempts to improve the properties of
linear subdivision schemes, have lead to nonlinear subdivision
schemes. For such schemes, the subdivision rules become data
dependant; in addition to the previously defined criteria, one
should add a stability property that ensures that the nonlinear
scheme is linearly affected by perturbations of the data (for
linear schemes, the stability is a direct consequence of the
convergence).

 For nonlinear subdivision  schemes, very few results concerning convergence or
stability are available, see for instance \cite{AL}, \cite{CDM},
\cite{DY}, \cite{Os}, \cite{DRS} and \cite{FM}.

A large family   of nonlinear subdivision schemes to which belong
the ENO, WENO or PPH schemes \cite{CDM},
\cite{ADLT} is made by the  schemes constructed as a perturbation of
 the four-point linear interpolatory Lagrange scheme based on centered
 degree $3$ polynomial interpolation. These schemes   are
  interpolatory subdivision schemes (i.e. based on interpolation rules)
   and are  constructed to avoid  the Gibbs oscillations occurring classically for
    linear interpolatory schemes (see Figure \ref{f1}).
    The schemes of this family  are unfortunately charac\-terized by a low regularity
    of the limit functions of type  $C^{1-}$. Moreover, the ENO scheme is unstable.

In \cite{DFH05}, a new linear
 four-point subdivision scheme was presented.
 Its refinement rule is  based on local cubic interpolation  followed
by a shift of $1/4$ or, in other words an evaluation at positions
$1/4$ and $3/4$ rather then the standard evaluation at $1/2$. This
new scheme was shown to be convergent towards a $C^2$ curve.

The aim of this paper is to analyze the scheme obtained using the
same trick (shift of $1/4$) for the PPH-type schemes  \cite{ADLT}
which are derived modifying the classical four point interpolatory
subdivision scheme substituting  the harmonic mean to the
arithmetic mean. After the definition of the new scheme in section
\ref{sec2} we analyze successively its convergence (section
\ref{sec3}), its  stability (section \ref{sec4}) and its order of
approximation (in section \ref{sec5}). Its  behavior in presence
of strongly variating data (Gibbs oscillations) is analyzed in
section  \ref{secG}.
  The last section is devoted to concluding remarks.

\section{A new nonlinear subdivision scheme}\label{sec2}

As mentioned above, the starting point of our work is the
construction of N. Dyn, M.S. Floater and K. Hormann in
\cite{DFH05}. There,   a new  linear four-point subdivision scheme
that generates $C^2$ curves is presented. Its refinement rule is
based on the local cubic Lagrange interpolation, followed by
evaluation at positions  $1/4$ and $3/4$ of the refined interval.
For all  $f\in l^\infty(\mathbb{Z})$, the scheme is then  given by
\begin{eqnarray}
\label{schemaDFH}
(Sf)_{2n} &= &   -\frac{7}{128}f_{n-1}+\frac{105}{128}f_{n}+\frac{35}{128}f_{n+1}-\frac{5}{128}f_{n+2},\nonumber\\
 (Sf)_{2n+1}& = &
-\frac{5}{128}f_{n-1}+\frac{35}{128}f_{n}+\frac{105}{128}f_{n+1}-\frac{7}{128}f_{n+2}.
\end{eqnarray}
Following \cite{ADLT} where a nonlinear scheme is derived
modifying the classical four-point interpolatory subdivision
scheme substituting  the harmonic mean to  the arithmetic mean, we
first obtain two new formulations of the scheme (\ref{schemaDFH}).
\begin{description}
\item[1]
\begin{eqnarray*}
(S f)_{2n} &= &
\frac{49}{64}f_{n}+\frac{14}{64}f_{n+1}+\frac{1}{64}f_{n+2}-\frac{7}{64}\frac{(d^2f_{n}+d^2f_{n+1})}{2},\\
 (S f)_{2n+1}& =
&\frac{15}{64}f_{n}+\frac{50}{64}f_{n+1}-\frac{1}{64}f_{n+2}-
\frac{5}{64}\frac{(d^2f_{n}+d^2f_{n+1})}{2}.
\end{eqnarray*}
\item[2]
\begin{eqnarray*}
(S f)_{2n} &=
&-\frac{1}{64}f_{n-1}+\frac{15}{64}f_{n}+\frac{50}{64}f_{n+1}-\frac{5}{64}\frac{(d^2f_{n}+d^2f_{n+1})}{2},
\\
(S f)_{2n+1}& =
&\frac{1}{64}f_{n-1}+\frac{49}{64}f_{n}+\frac{14}{64}f_{n+1}-\frac{7}{64}\frac{(d^2f_{n}+d^2f_{n+1})}{2}.
\end{eqnarray*}
\end{description}
where $(d^2f)$ is defined by $d^2f_n=f_{n+1}-2f_n+f_{n-1}$.

The two formulations differs essentially in the distribution of
the points $f_n$ contributing to the three first terms of {\bf 1}
and {\bf 2}.

Using the same strategy as in \cite{ADLT},  we define the new
nonlinear subdivision scheme $S_{\pphapprox}$ associated to
(\ref{schemaDFH}) by
\begin{description}
\item[If $|d^2f_n|\geq |d^2f_{n+1}|$,]
\begin{eqnarray*} \label{schemaPPHAPPROX} (S_\pphapprox f)_{2n} &= &
\frac{49}{64}f_{n}+\frac{14}{64}f_{n+1}+\frac{1}{64}f_{n+2}-\frac{7}{64}\pph(d^2f_{n},d^2f_{n+1}),\\
\nonumber (S_\pphapprox f)_{2n+1}& =
&\frac{15}{64}f_{n}+\frac{50}{64}f_{n+1}-\frac{1}{64}f_{n+2}-\frac{5}{64}\pph(d^2f_{n},d^2f_{n+1}).
\end{eqnarray*}
\item[If $|d^2f_n|<|d^2f_{n+1}|$,]
\begin{eqnarray*}
(S_\pphapprox f)_{2n} &=
&-\frac{1}{64}f_{n-1}+\frac{15}{64}f_{n}+\frac{50}{64}f_{n+1}-\frac{5}{64}\pph(d^2f_{n},d^2f_{n+1}),
\\
(S_\pphapprox f)_{2n+1}& =
&\frac{1}{64}f_{n-1}+\frac{49}{64}f_{n}+\frac{14}{64}f_{n+1}-\frac{7}{64}\pph(d^2f_{n},d^2f_{n+1}),
\end{eqnarray*}
\end{description}
where $\pph$ stands for  the harmonic mean defined by
\begin{equation}\label{harmonic}
\nonumber (x,y) \in I\!\!R^2 \mapsto {\pph}(x,y):= \frac{xy}{x+y} ({\rm sgn}(xy) +1),
\end{equation}
with ${\rm sgn}(x)=1$ if $x \geq 0$ and ${\rm sgn}(x)=-1$ if
$x<0$.

The initial motivation for the substitution of the arithmetic mean by the harmonic mean is the elimination of oscillations near strong variating data thanks to the fact that

\begin{equation}\label{meandos}
|\pph(x,y)| \leq 2 \min(|x|, |y|),
\end{equation}
substitutes to
\begin{eqnarray*}
\frac{x+y}{2}\leq max(|x|,|y|).
\end{eqnarray*}


Before analyzing in details the properties of the new scheme $S_{\pphapprox}$  we summarize the most important properties of the harmonic mean in
the following proposition (see \cite{AL} for more details).

\begin{proposition}
\label{prop2.1} For all $(x,y)\in \mathbb{R}^2$, the harmonic mean
$\pph(x,y) $ satisfies
\begin{enumerate}
\item $\pph(x,y)=\pph(y,x)$.
 \item $\pph(x,y)=0 \quad if\; xy\leq0$.
\item $\pph(-x,-y)=-\pph(x,y)$. \item $\pph(x,y)  =
\frac{sign(x)+sign(y)}{2} min(|x|,|y|)\left [ 1+\left |
\frac{x-y}{x+y} \right | \right].$ \item $|\pph(x,y)|\leq
\max{(|x|,|y|)}$. \item $|\pph(x,y)| \leq 2 \min{(|x|,|y|)}$.
\item For $x,y>0$, $\min(x,y) \leq \pph(x,y)\leq \frac{x+y}{2}$.
\item If $x=O(1)$, $y=O(1)$, $|y-x|=O(h)$ and $xy>0$ then
$$|\frac{x+y}{2}-\pph(x,y)|=O(h^2).$$ \item
$|\pph(x_1,y_1)-\pph(x_2,y_2)|\leq 2 \max (|x_1-x_2|,|y_1-y_2|)$.
\end{enumerate}
\end{proposition}

\section{Convergence and Regularity}\label{sec3}

We recall the following definition.

\begin{definition}
\label{defconv}
 A subdivision scheme $S$  is said to be convergent  if
\begin{equation}
 \forall f \in l^\infty(\mathbb{Z}),\, \exists
 S^\infty f \in C^0(\mathbb{R}) \mbox{ such that }
\lim_{j \rightarrow +\infty}  sup_{n\in \mathbb{Z}}|(S^jf)_n-S^\infty f (n2^{-j})|=0.
\end{equation}
\end{definition}

In order to derive the convergence, we rewrite  the nonlinear
subdivision scheme $S_\pphapprox$ as a perturbation of a classical
two-point linear subdivision scheme, $S_\chaikin$, introduced by
G. Chaikin in \cite{chaikin} and defined by

\begin{eqnarray}
\label{schemachaikin}
(S_\chaikin f)_{2n} &= &
\frac{3}{4}f_{n}+\frac{1}{4}f_{n+1},\\
\nonumber
(S_\chaikin f)_{2n+1}& =
&\frac{1}{4}f_{n}+\frac{3}{4}f_{n+1}.
\end{eqnarray}
The scheme $S_\chaikin$ is known to be convergent  with a regularity $C^{2-}$.

Writing
\begin{description}
\item[If $|d^2f_n|\geq|d^2f_{n+1}|$,]
\begin{eqnarray*}
\label{schemaDFHnlpertub} (S_\pphapprox f)_{2n} &= &
\frac{3}{4}f_{n}+\frac{1}{4}f_{n+1}+\frac{1}{64}d^2f_{n+1}-\frac{7}{64}\pph(d^2f_{n},d^2f_{n+1}),\\
\nonumber (S_\pphapprox f)_{2n+1}& =
&\frac{1}{4}f_{n}+\frac{3}{4}f_{n+1}-\frac{1}{64}d^2f_{n+1}-\frac{5}{64}\pph(d^2f_{n},d^2f_{n+1}),
\end{eqnarray*}
\item[If $|d^2f_n|<|d^2f_{n+1}|$,]
\begin{eqnarray*}
(S_\pphapprox f)_{2n} &=
&\frac{3}{4}f_{n}+\frac{1}{4}f_{n+1}-\frac{1}{64}d^2f_{n}-\frac{5}{64}\pph(d^2f_{n},d^2f_{n+1}),
\\
(S_\pphapprox f)_{2n+1}& =
&\frac{1}{4}f_{n}+\frac{3}{4}f_{n+1}+\frac{1}{64}d^2f_{n}-\frac{7}{64}\pph(d^2f_{n},d^2f_{n+1}),
\end{eqnarray*}
\end{description}

we get that  $S_{\pphapprox}$ can be expressed as
$$
S_\pphapprox f=S_\chaikin f +F(d^2f),
$$
with
\begin{equation}
\label{Fpphapprox1} F(d^2f)_{2n}= \left \{
\begin{array}{ll}
\frac{1}{64}d^2f_{n+1}-\frac{7}{64}\pph(d^2f_{n},d^2f_{n+1}) & \quad \textrm{ if } |d^2f_n|>|d^2f_{n+1}|,\\
-\frac{1}{64}d^2f_{n+1}-\frac{5}{64}\pph(d^2f_{n},d^2f_{n+1}) &
\quad \textrm{ if } |d^2f_n|<|d^2f_{n+1}|,
\end{array}
\right .
\end{equation}
and
\begin{equation}
\label{Fpphapprox2} F(d^2f)_{2n+1}= \left \{
\begin{array}{ll}
-\frac{1}{64}d^2f_{n+1}-\frac{5}{64}\pph(d^2f_{n},d^2f_{n+1}) & \quad \textrm{ if } |d^2f_n|>|d^2f_{n+1}|,\\
\frac{1}{64}d^2f_{n+1}-\frac{7}{64}\pph(d^2f_{n},d^2f_{n}) & \quad
\textrm{ if } |d^2f_n|<|d^2f_{n+1}|.
\end{array}
\right .\\
\end{equation}
\ \\
\ \\
To analyze the convergence of $S_\pphapprox$, we use a result
proved in  \cite{ADL06}, \cite{ADL} that reads:

A sufficient condition  for the convergence of a nonlinear
subdivision scheme $S_{NL}: l^{\infty}(\mathbb{Z}) \rightarrow
l^{\infty}(\mathbb{Z})$ of the form:
\begin{equation}\label{perturbation}
\forall f \in l^{\infty}(\mathbb{Z}), \quad \forall n\in
\mathbb{Z} \qquad \left \lbrace \begin{array}{lll}
(S_{NL}f)_{2n+1} $=$ (Sf)_{2n+1}+F(\delta f)_{2n+1},\\
(S_{NL}f)_{2n} $=$ f_{n},
\end{array} \right .
\end{equation}
where $F$ is a nonlinear operator defined on
$l^{\infty}({\mathbb{Z}})$, $\delta$ is a  linear  and continuous
operator
 on $l^{\infty}(\mathbb{Z})$ and $S$ is a linear and convergent subdivision scheme is:

\begin{theorem}
 \label{th2.1}

  If $F, S $ and $\delta$ given in (\ref{perturbation}) verify:

\begin{eqnarray}
 & & \exists M>0 \quad  \textrm{such that} \quad  \forall d \in
 l^{\infty}({\mathbb{Z}})
\quad \; || F(d)||_\infty \leq M ||d||_{\infty}, \label{h1}\\
 & & \exists c<1 \; \textrm{such that} \qquad
|| \delta S(f)+ \delta F(\delta f)  ||_\infty \leq c || \delta
f||_\infty, \label{h2}
\end{eqnarray}
then the subdivision  scheme $S_{NL}$ is uniformly convergent.
Moreover, if $S$ is $C^{\alpha^-}$ convergent
 then, for all sequence $f \in l^{\infty} (\mathbb{Z}),\;S^{\infty}_{NL}(f) \in C^{\beta^-}$
 with $\beta=\min{\left (\alpha,-log_2(c) \right )}$. \\
\end{theorem}




Using theorem \ref{th2.1}, we are going to prove the following
result.

\begin{theorem}
\label{th2.2} The nonlinear subdivision scheme $S_{\pphapprox}$ is convergent
with a regularity at least
 $C^{\beta-}$
with $\beta\geq -log_2(\frac{7}{16})>1$.
\end{theorem}

{\bf Proof}

 From the properties  of the harmonic mean (Proposition  \ref{prop2.1}),
 \begin{equation}
 \label{eq1demoth2.7}
 |c_1d^2f_{n+1}-c_2\pph(d^2f_{n},d^2f_{n+1})| \leq \max{(c_1,c_2)} ||d^2f||_\infty.
 \end{equation}
 \ \\
For the perturbation $F$ defined in   (\ref{Fpphapprox1}) and
(\ref{Fpphapprox2}), it is then easy to  see that
for all $d \in l^\infty (\mathbb{Z})$,
\begin{equation}
\label{eq0demoth2.7} ||F(d)||_\infty \leq
\frac{7}{64}||d||_\infty,
\end{equation}
that is hypothesis  (\ref{h1}).
\ \\
We now consider hypothesis $(\ref{h2})$ related, in this case,  to
the contraction of  the second
order differences $(d^2f)$. To simplify the notations we call $f^1=S_\pphapprox(f)$.\\
Different cases must be considered:
\begin{description}
\item[Case 1: k=2n+1,] study of
$(d^2f^1)_{2n+1}=f^1_{2n+2}-2f^1_{2n+1}+f^1_{2n}$
\begin{center}
\begin{tabular}{ll}
case 1A$_1$:&  $|d^2f_{n}|\geq|d^2f_{n+1}|$ and $|d^2f_{n+1}|\geq|d^2f_{n+2}|$,\\
case 1A$_2$:&  $|d^2f_{n}|<|d^2f_{n+1}|$ and $|d^2f_{n+1}|<|d^2f_{n+2}|$,\\
case 1B$_1$:&  $|d^2f_{n}|\geq|d^2f_{n+1}|$ and $|d^2f_{n+1}|<|d^2f_{n+2}|$,\\
case 1B$_2$:&  $|d^2f_{n}|<|d^2f_{n+1}|$ and
$|d^2f_{n+1}|\geq|d^2f_{n+2}|$.
\end{tabular}
\end{center}
\item[Case 2: k=2n,] study of
$(d^2f^1)_{2n}=f^1_{2n+1}-2f^1_{2n}+f^1_{2n-1}$
\begin{center}
\begin{tabular}{ll}
case 2A$_1$:&  $|d^2f_{n}|\geq|d^2f_{n+1}|$ and $|d^2f_{n-1}|\geq|d^2f_{n}|$,\\
case 2A$_2$:&  $|d^2f_{n}|<|d^2f_{n+1}|$ and $|d^2f_{n-1}|<|d^2f_{n}|$,\\
case 2B$_1$:&  $|d^2f_{n}|\geq|d^2f_{n+1}|$ and $|d^2f_{n-1}|<|d^2f_{n}|$,\\
case 2B$_2$:&  $|d^2f_{n}|<|d^2f_{n+1}|$ and
$|d^2f_{n-1}|\geq|d^2f_{n}|$.
\end{tabular}
\end{center}
\end{description}
The others cases follow by symmetry.\\
\ \\

$\bullet$ Cases 1A: We obtain for the case 1A$_1$
\begin{eqnarray*}
\nonumber (d^2f^1)_{2n+1} & = & \frac{1}{4}f_{n+2}-
\frac{2}{4}f_{n+1}+ \frac{1}{4}f_{n} +
\frac{1}{64}d^2f_{n+2}+\frac{3}{64}d^2f_{n+1}\\
\nonumber
 & & -\frac{7}{64}\pph(d^2f_{n+1},d^2f_{n+2})+\frac{3}{64}\pph(d^2f_{n},d^2f_{n+1})\\
 \nonumber
 &=& \frac{1}{4}d^2f_{n+1}+
\frac{1}{64}d^2f_{n+2}+\frac{3}{64}d^2f_{n+1}\\ \nonumber
& &-\frac{7}{64}\pph(d^2f_{n+1},d^2f_{n+2})+\frac{3}{64}\pph(d^2f_{n},d^2f_{n+1})\\
 \nonumber &=& \frac{19}{64}d^2f_{n+1}+
\frac{1}{64}d^2f_{n+2}   \nonumber \\ && \label{eq1A1demoth2.7}
-\frac{7}{64}\pph(d^2f_{n+1},d^2f_{n+2})+\frac{3}{64}\pph(d^2f_{n},d^2f_{n+1}).
\end{eqnarray*}
 Using  equation (\ref{eq1demoth2.7}) for $d^2f_{n+1}$ and
  $\pph(d^2f_{n+1},d^2f_{n+2})$ and Proposition  \ref{prop2.1},
  we have
 \begin{eqnarray}
\label{eq2demoth2.7}
 |(d^2f^1)_{2n+1}| & \leq & \frac{19+1+3}{64}||d^2f||_\infty
   \leq   \frac{23}{64}||d^2f||_\infty.
 \end{eqnarray}
 Similarly for the case 1A$_2$, we have
 \begin{eqnarray*}
(d^2f^1)_{2n+1} & = & \frac{15}{64}d^2f_{n+1}-
\frac{3}{64}d^2f_{n}-\frac{5}{64}\pph(d^2f_{n+1},d^2f_{n+2})\nonumber \\
 \label{eq1A2demoth2.7} & &+\frac{9}{64}\pph(d^2f_{n},d^2f_{n+1}),
\end{eqnarray*}
and (\ref{eq2demoth2.7}) remains valid.

 \ \\

$\bullet$ Cases 1B: We obtain for the case 1B$_1$
\begin{eqnarray*}
\nonumber (d^2f^1)_{2n+1} & = & \frac{1}{4}f_{n+2}-
\frac{2}{4}f_{n+1}+ \frac{1}{4}f_{n} -
\frac{1}{64}d^2f_{n+1}+\frac{3}{64}d^2f_{n+1}\\
\nonumber
 & & -\frac{5}{64}\pph(d^2f_{n+1},d^2f_{n+2})+\frac{3}{64}\pph(d^2f_{n},d^2f_{n+1})\\
   \nonumber
  &=& \frac{18}{64}d^2f_{n+1}\\
  \label{eq1B1demoth2.7}
&&-\frac{5}{64}\pph(d^2f_{n+1},d^2f_{n+2})+\frac{3}{64}\pph(d^2f_{n},d^2f_{n+1}).
  \end{eqnarray*}

Using  equation (\ref{eq1demoth2.7}) for $d^2f_{n+1}$ and
$\pph(d^2f_{n+1},d^2f_{n+2})$, and   Proposition
 \ref{prop2.1}, we have
 \begin{eqnarray}
 \label{eq3demoth2.7}
 |(d^2f^1)_{2n+1}| & \leq & \frac{18+3}{64}||d^2f||_\infty
\leq   \frac{21}{64}||d^2f||_\infty.
 \end{eqnarray}
 Similarly for the case 1B$_2$, we have
 \begin{eqnarray*}
   (d^2f^1)_{2n+1}
  &=& \frac{16}{64}d^2f_{n+1}+\frac{1}{64}d^2f_{n+2}
  -\frac{3}{64}d^2f_{n}\nonumber \\ \label{eq1B2demoth2.7} &&-\frac{7}{64}\pph(d^2f_{n+1},d^2f_{n+2})
  +\frac{9}{64}\pph(d^2f_{n},d^2f_{n+1}),
  \end{eqnarray*}
and (\ref{eq3demoth2.7}) remains valid.
\ \\

$\bullet$ Cases 2A: We obtain for the case 2A$_1$
\begin{eqnarray*}
\nonumber (d^2f^1)_{2n} & = & \frac{1}{4}f_{n+1}-
\frac{2}{4}f_{n}+ \frac{1}{4}f_{n-1} -
\frac{3}{64}d^2f_{n+1}-\frac{1}{64}d^2f_{n}\\
\nonumber
 & & +\frac{9}{64}\pph(d^2f_{n},d^2f_{n+1})-\frac{5}{64}\pph(d^2f_{n-1},d^2f_{n})\\
 \nonumber
 &=& \frac{1}{4}d^2f_{n}-
\frac{3}{64}d^2f_{n+1}-\frac{1}{64}d^2f_{n}\\  \nonumber & &
+\frac{9}{64}\pph(d^2f_{n},d^2f_{n+1})-\frac{5}{64}\pph(d^2f_{n-1},d^2f_{n})\\
 \nonumber
 &=& \frac{15}{64}d^2f_{n}-
\frac{3}{64}d^2f_{n+1}\\
&& \label{eq2A1demoth2.7}
 +\frac{9}{64}\pph(d^2f_{n},d^2f_{n+1})
-\frac{5}{64}\pph(d^2f_{n-1},d^2f_{n}).
\end{eqnarray*}
Using equation (\ref{eq1demoth2.7}) for $d^2f_{n+1}$ and
$\pph(d^2f_{n},d^2f_{n+1})$, and for $d^2f_{n-1}$ and
$\pph(d^2f_{n-1},d^2f_{n})$, we have
 \begin{eqnarray}
 \label{eq4demoth2.7}
 |(d^2f^1)_{2n}| & \leq & \frac{15+9}{64}||d^2f||_\infty
   \leq   \frac{28}{64}||d^2f||_\infty.
 \end{eqnarray}
 Similarly for the case 2A$_2$, we have
  \begin{eqnarray*}
(d^2f^1)_{2n} & = &  \frac{19}{64}d^2f_{n}+ \frac{1}{64}d^2f_{n-1}
\nonumber
\\
\label{eq2A2demoth2.7}
&&+\frac{3}{64}\pph(d^2f_{n},d^2f_{n+1})-\frac{7}{64}\pph(d^2f_{n-1},d^2f_{n}),
\end{eqnarray*}
and (\ref{eq4demoth2.7}) remains valid.
 \ \\

 $\bullet$ Cases 2B: We obtain for the case 2B$_1$
 \begin{eqnarray*}
 \nonumber
(d^2f^1)_{2n} & = & \frac{1}{4}f_{n+1}- \frac{2}{4}f_{n}+
\frac{1}{4}f_{n-1} -
\frac{3}{64}d^2f_{n+1}+\frac{1}{64}d^2f_{n-1}\\
\nonumber
 & & +\frac{9}{64}\pph(d^2f_{n},d^2f_{n+1})-\frac{7}{64}\pph(d^2f_{n-1},d^2f_{n})\\
\nonumber
 &= &  \frac{1}{4}d^2f_{n}-
\frac{3}{64}d^2f_{n+1}+\frac{1}{64}d^2f_{n-1}\\
\label{eq2B1demoth2.7} &&
+\frac{9}{64}\pph(d^2f_{n},d^2f_{n+1})-\frac{7}{64}\pph(d^2f_{n-1},d^2f_{n}).
\end{eqnarray*}
Using equation (\ref{eq1demoth2.7}) for $d^2f_{n+1}$ and
$\pph(d^2f_{n},d^2f_{n+1})$, and for
 $d^2f_{n}$ and $\pph(d^2f_{n-1},d^2f_{n})$, we have
 \begin{eqnarray}
 \label{eq5demoth2.7}
 |(d^2f^1)_{2n}| & \leq & \frac{16+9+1}{64}||d^2f||_\infty
 \leq   \frac{28}{64}||d^2f||_\infty.
 \end{eqnarray}
  Similarly for the case 2B$_2$, we have
\begin{eqnarray*}
\label{eq2B2demoth2.7} (d^2f^1)_{2n} & = & \frac{18}{64}d^2f_{n}
\nonumber  \\ && +
\frac{3}{64}\pph(d^2f_{n},d^2f_{n+1})-\frac{5}{64}\pph(d^2f_{n-1},d^2f_{n}),
\end{eqnarray*}
and (\ref{eq5demoth2.7}) remains valid.

 From  equations (\ref{eq2demoth2.7}), (\ref{eq3demoth2.7}), (\ref{eq4demoth2.7})
 and
 (\ref{eq5demoth2.7}), we deduce that for all $f\in l^\infty(\mathbb{Z})$
 \begin{equation}
 \label{eq6demoth2.7}
 ||d^2S_{\pphapprox}f||_\infty \leq \frac{7}{16}||d^2f||_\infty.
 \end{equation}

  Therefore $S_{\pphapprox}$ verifies the hypothesis (\ref{h2}) of the Theorem
   \ref{th2.1}. In particular, we obtain the convergence of $S_{\pphapprox}$.\\

   For the regularity, we use again  Theorem \ref{th2.1}. According to the values $\alpha=2$ and $c=\frac{7}{16}$ we  obtain the regularity
   constant  $\beta=\min{(2,-log_2\left ( \frac{7}{16})
\right )}\approx 1.192$.

$\Box$

{\bf Numerical Regularity}

Following \cite{Kuijt98},  the regularity of a limit function can be evaluated  numerically.
Using $S_1$ and $S_2$  the subdivision schemes for the differences of order $1$ and
$2$ associated to $S_\pphapprox$ (that can be derived due to the specific definition of  $S_\pphapprox$), the following quantities are  estimated for $k=1,2$,
  $$-log_2\left(2^k\frac{||(S^{j+1}_kf)_{n+1}-(S^{j+1}_kf)_n ||_\infty}
  {||(S^{j}_kf)_{n+1}-(S^{j}_kf)_n
||_\infty} \right ).$$ They provide an estimate for $\beta_1$ and $\beta_2$ such that the limit function belongs to $C^{1+\beta_1-}$ and $C^{2+\beta_2-}$. From table $\ref{table4ch2}$,
the numerical estimate of the regularity is  $C^{2.438-}$. We recall that  the corresponding estimate for
the  linear scheme \cite{DFH05}  is
  $C^{2.67-}$. \\
    \begin{table}[!h]
 \centering
\begin{tabular}{|c||c|c|c|c|c|c|}
\hline
 $j$ &  5 & 6&  7   & 8  & 9 & 10 \\
 \hline
 $\beta_1$ & 0.9999   &  0.9999  &  1   &  1 &  1 & 1\\
\hline
$\beta_2$ &   0.4395  &    0.7738  &   1.2615 &   0.6541  &   0.4387 &   0.4388 \\
\hline
\end{tabular}
\caption{ Numerical estimates of the limit function regularity
$C^{1+\beta_1-}$ and $C^{2+\beta_2-}$ for $S_\pphapprox$. }
 \label{table4ch2}
\end{table}

\section{Stability}\label{sec4}
For simplicity in notations we call for any initial sequence $f^0$, and any
 $j \in \mathbb{N}$, $f^{j+1}=S(f^j)$. We recall the following
 definition.

\begin{definition}
\label{def1.9}
 A convergent subdivision scheme is stable if
\begin{equation}
 \exists C <+\infty  \textrm{ such that } \forall f^0,g^0\in
 l^{\infty}(\mathbb{Z})
 \quad ||S^\infty f -S^\infty g ||_{L^\infty} \leq C
||f^0-g^0||_{l^\infty}.
\end{equation}
\end{definition}

As for the convergence, to derive the stability of $S_\pphapprox$
we use the following theorem of \cite{ADL}.

\begin{theorem}\label{th2.3}

If $F, S$ and $\delta$ given in (\ref{perturbation}) verify:
$\exists M>0, c<1$ such that $\forall f,g, d_1,d_2$,
\begin{eqnarray}
  \quad  ||F(d_1)-F(d_2)||_\infty  \leq M ||d_1-d_2||_{\infty}, \label{h1th2.3}\\
  \quad  \Vert \delta (S_{NL}f-S_{NL}g)  \Vert_\infty \leq c \Vert \delta (f-g)
\Vert_\infty, \label{h2th2.3}
\end{eqnarray}
then  the  nonlinear subdivision scheme $S_{NL}$ is stable.
\end{theorem}

In order to check the hypotheses of  Theorem \ref{th2.3} for
$S_{NL}=S_\pphapprox$ we first prove the following lemma.

\begin{lemma}
\label{lemme2.1} Let be $f,\, g \in l^\infty(\mathbb{Z})$, if
$|d^2f_n|\geq|d^2f_{n+1}|$ and $|d^2g_{n+1}|\geq|d^2g_n|$ then
$$
|d^2f_{n+1}+d^2g_n-2\pph(d^2g_n,d^2g_{n+1})|\leq
3||d^2f-d^2g||_\infty.
$$
\end{lemma}

{\bf Proof}

We consider different cases.\\
\ \\
$\bullet$ If $d^2g_{n+1}d^2g_n<0$, we have $$
d^2f_{n+1}+d^2g_n-2\pph(d^2g_n,d^2g_{n+1})=d^2f_{n+1}+d^2g_n.
$$
\begin{itemize}
\item[$\ast$] if $d^2f_{n+1}d^2g_n<0$, using that
$|d^2f_n|\geq|d^2f_{n+1}|$
$$
|d^2f_{n+1}+d^2g_n|\leq |d^2f_{n+1}-d^2g_n| \leq |d^2f_n-d^2g_n|.
$$
\item[$\ast$] if $d^2f_{n+1}d^2g_n\geq0$, we have
$d^2f_{n+1}d^2g_{n+1}<0$ and with $|d^2g_{n+1}|\geq|d^2g_n|$, we
obtain
$$
|d^2f_{n+1}+d^2g_n|\leq |d^2f_{n+1}-d^2g_{n+1}|.
$$
\end{itemize}
\ \\
$\bullet$ If $d^2g_{n+1}d^2g_n\geq0$, we recall (Proposition
\ref{prop2.1}) that if
$x,y>0$, $\min{(x,y)}\leq \pph(x,y) \leq max{(x,y)}$.\\
Without loss of generality, we suppose that $d^2g_n\geq0$. We
denote
$$H=d^2f_{n+1}+d^2g_n-2\pph(d^2g_n,d^2g_{n+1}).$$
\begin{itemize}
\item[$\ast$] if $H>0$,
\begin{eqnarray*}
H & \leq &  d^2f_{n+1}+d^2g_n-2\min{(|d^2g_n|,|d^2g_{n+1}|)}\\
 & \leq & d^2f_{n+1}-d^2g_n\\
 & \leq &|d^2f_n-d^2g_n|.\\
\end{eqnarray*}
\item[$\ast$] if $H<0$,
\begin{eqnarray*}
H & \geq & d^2f_{n+1}+d^2g_n-2\max{(|d^2g_n|,|d^2g_{n+1}|)}    \\
 & \geq & d^2f_{n+1}+d^2g_n-2d^2g_{n+1}\\
  & \geq & d^2f_{n+1}-d^2g_{n+1}+d^2g_n-d^2g_{n+1}.
  \end{eqnarray*}
 We have again to consider different cases according to the sign of $d^2f_n$.
  \begin{itemize}
  \item $d^2f_n\geq0$, we have $d^2f_n-d^2f_{n+1}\geq0$
  \begin{eqnarray*}
H &\geq &(d^2f_{n+1}-d^2g_{n+1})+(d^2g_n-d^2g_{n+1})-(d^2f_n-d^2f_{n+1}) \\
 &\geq & 2(d^2f_{n+1}-d^2g_{n+1})+(d^2g_n-d^2f_{n})\geq -3||d^2f-d^2g||_\infty.
 \end{eqnarray*}
 \item  $d^2f_n<0$ and $d^2f_{n+1}d^2f_n\geq0$, then $d^2f_{n+1}<0$
   \begin{eqnarray*}
H &\geq & (d^2f_{n+1}-d^2g_{n+1})+(d^2f_{n+1}-d^2g_{n+1}) \\
 &\geq &2(d^2f_{n+1}-d^2g_{n+1})\geq -2||d^2f-d^2g||_\infty.
 \end{eqnarray*}
 \item  $d^2f_n<0$ and $d^2f_{n+1}d^2f_n<0$, we have from hypothesis that
 $d^2f_{n+1}+d^2f_n<0$,
\begin{eqnarray*}
H &\geq &(d^2f_{n+1}-d^2g_{n+1})+(d^g_{n}-d^2g_{n+1})+d^2f_{n+1}+d^2f_{n} \\
 &\geq &2(d^2f_{n+1}-d^2g_{n+1})+d^2f_n+d^2g_n\\
 &\geq &2(d^2f_{n+1}-d^2g_{n+1})+d^2f_n-d^2g_n \geq  -3||d^2f-d^2g||_\infty.
 \end{eqnarray*}

$\Box$

\end{itemize}
\end{itemize}

We are now ready to prove the stability of $S_\pphapprox$.

\begin{theorem}
\label{th2.8} The scheme $S_\pphapprox$ is stable.
\end{theorem}

{\bf Proof}

We check  the hypotheses of  Theorem \ref{th2.3}.

Firstly, we start with the hypothesis (\ref{h1th2.3}) for $F$.\\
 Using the expressions of perturbation $F$,
(\ref{Fpphapprox1}) and (\ref{Fpphapprox2}), and  Proposition
 \ref{prop2.1}, we obtain for all $d_1,d_2 \in
l^\infty(\mathbb{Z}) $ that
\begin{eqnarray*}
||F(d_1)-F(d_2)||_\infty  & \leq  &\frac{1+7\cdot 2}{64} ||d_1-d_2||_\infty.\\
\end{eqnarray*}

Secondly, we have to verify the contraction hypothesis
(\ref{h2th2.3}).

 For a couple  $f$, $g\in l^\infty(\mathbb{Z})$,
we
 study $(d^2f^1-d^2g^1)_{k}$ for $k=2n+1$ (case 1) or $k=2n$ (case 2).

 We consider different cases, according to the  proof
of  Theorem \ref{th2.2} for  $f$ or $g$.
\ \\
\textbf{For k=2n+1}, we have $7$ cases to study, see table \ref{table5ch2}, the others
cases being  deduced by symmetry.\\
   \begin{table}[!h]
 \centering
\begin{tabular}{|c|c|c|c||c|c|c|c|}
\hline
 & $f$ verifies & $g$ verifies  & notation &  & $f$ verifies  & $g$ verifies  & notation \\
 \hline
 \hline
\multirow{4}{1cm}{Case 1A} &  1A$_1$ &  1A$_1$ & $1A_1-1A_1$ &  \multirow{4}{1cm}{Case 1B}   &  1B$_1$  &  1A$_2$ &$1B_1-1A_2$\\
\cline{2-4}\cline{6-8}
 &  1A$_1$ & 1A$_2$   & $1A_1-1A_2$ &   &1B$_1$  &  1B$_1$ &$1B_1-1B_1$\\
\cline{2-4}\cline{6-8}
 &  1A$_1$ &  1B$_1$ & $1A_1-1B_1$&  & 1B$_1$  &  1B$_2$ &$1B_1-1B_2$\\
\cline{2-4}\cline{6-8}
 &  1A$_1$ &   1B$_2$  &    $1A_1-1B_2$  &    & 1B$_2$  &  1A$_2$ & $1B_2-1A_2$  \\
\hline
\end{tabular}
\caption{ Cases to consider for $k=2n+1$ in the proof of the
stability of $S_\pphapprox$. } \label{table5ch2}
\end{table}

\ \\
$\bullet$ Case $1A_1-1A_1$: with  equation
(\ref{eq1A1demoth2.7}) we can obtain directly
\begin{eqnarray*}
|d^2f^1_{2n+1}-d^2g^1_{2n+1}| & \leq & \frac{19+1+7\cdot2+3\cdot2}{64}
||d^2f-d^2g||_\infty\\
 & \leq &   \frac{5}{8}||d^2f-d^2g||_\infty.\\
\end{eqnarray*}

$\bullet$ Case $1A_1-1A_2$: from  equations
(\ref{eq1A1demoth2.7}) and (\ref{eq1A2demoth2.7}), we obtain
\begin{eqnarray*}
d^2f^1_{2n+1}-d^2g^1_{2n+1} &= & \frac{19}{64}d^2f_{n+1}-
\frac{15}{64}d^2g_{n+1}+
\frac{1}{64}d^2f_{n+2}+ \frac{3}{64}d^2g_{n}-\frac{7}{64}\pph(d^2f_{n+1},d^2f_{n+2})\\
 & &
+\frac{5}{64}\pph(d^2g_{n+1},d^2g_{n+2})
+\frac{3}{64}\pph(d^2f_{n},d^2f_{n+1}) +\frac{9}{64}\pph(d^2g_{n},d^2g_{n+1})\\
 & =& \frac{16}{64}(d^2f_{n+1}-d^2g_{n+1})-\frac{7}{64}(\pph(d^2f_{n+1},d^2f_{n+2})-\pph(d^2g_{n+1},d^2g_{n+2}))\\
 & & +\frac{1}{64}(d^2f_{n+2}+d^2g_{n+1}-2\pph(d^2g_{n+1},d^2g_{n+2}))\\
 &&
 +\frac{3}{64}(\pph(d^2f_{n},d^2f_{n+1})-\pph(d^2g_{n+1},d^2g_{n+2}))\\
  & & +\frac{3}{64}(d^2g_{n}+d^2f_{n+1}-2\pph(d^2g_{n},d^2g_{n+1})) .
\end{eqnarray*}
Applying  Lemma \ref{lemme2.1} and  Proposition
\ref{prop2.1}, we obtain
\begin{eqnarray*}
|d^2f^1_{2n+1}-d^2g^1_{2n+1}| &\leq  &
\frac{16+7\cdot2+2+1+3\cdot2+3\cdot3}{64}
||d^2f-d^2g||_\infty\\
 & \leq & \frac{45}{64}||d^2f-d^2g||_\infty.\\
\end{eqnarray*}

Similarly

$\bullet$ Case $1A_1-1B_1$: from  equations (\ref{eq1A1demoth2.7})
and (\ref{eq1B1demoth2.7}), we obtain applying the Lemma
\ref{lemme2.1} and the Proposition \ref{prop2.1}
\begin{eqnarray*}
|d^2f^1_{2n+1}-d^2g^1_{2n+1}| & \leq & \frac{21}{32}||d^2f-d^2g||_\infty.\\
\end{eqnarray*}

$\bullet$ Case $1A_1-1B_2$: from  equations (\ref{eq1A1demoth2.7})
and (\ref{eq1B2demoth2.7}), we obtain applying the Lemma
\ref{lemme2.1} and the Proposition \ref{prop2.1}
\begin{eqnarray*}
|d^2f^1_{2n+1}-d^2g^1_{2n+1}| & \leq & \frac{23}{32}||d^2f-d^2g||_\infty.\\
\end{eqnarray*}

$\bullet$ Case $1B_1-1A_2$: from  equations (\ref{eq1B1demoth2.7})
and (\ref{eq1A2demoth2.7}), we obtain applying the Lemma
\ref{lemme2.1} and the Proposition \ref{prop2.1}
\begin{eqnarray*}
|d^2f^1_{2n+1}-d^2g^1_{2n+1}|& \leq & \frac{5}{8}||d^2f-d^2g||_\infty.\\
\end{eqnarray*}

$\bullet$ Case $1B_1-1B_1$: from  equation
(\ref{eq1B1demoth2.7}) we conclude
\begin{eqnarray*}
|d^2f^1_{2n+1}-d^2g^1_{2n+1}| & \leq & \frac{18+5\cdot2+3\cdot2}{64}||d^2f-d^2g||_\infty\\
 & \leq &   \frac{17}{32}||d^2f-d^2g||_\infty.\\
\end{eqnarray*}

$\bullet$ Case $1B_1-1B_2$: from  equations (\ref{eq1B1demoth2.7})
and (\ref{eq1B2demoth2.7}), we obtain applying the Lemma
\ref{lemme2.1} and the Proposition \ref{prop2.1}
\begin{eqnarray*}
|d^2f^1_{2n+1}-d^2g^1_{2n+1}| & \leq & \frac{3}{8}||d^2f-d^2g||_\infty.\\
\end{eqnarray*}

$\bullet$ Case $1B_2-1A_2$: from equations (\ref{eq1B2demoth2.7})
and (\ref{eq1A2demoth2.7}), we obtain applying the Lemma
\ref{lemme2.1} and the Proposition \ref{prop2.1}
\begin{eqnarray*}
|d^2f^1_{2n+1}-d^2g^1_{2n+1}| & \leq & \frac{11}{16}||d^2f-d^2g||_\infty.\\
\end{eqnarray*}
\ \\
\textbf{For k=2n}, we have $7$ other cases to study, see table
\ref{table6ch2}. The others cases are deduced by symmetry.\\
 We notice that in
 equations (\ref{eq1A2demoth2.7}) and (\ref{eq2A1demoth2.7})
the same  coefficients, but the cases are not completely
equivalent.\\
   \begin{table}[!h]
 \centering
\begin{tabular}{|c|c|c|c||c|c|c|c|}
\hline
 & $f$ verifies & $g$ verifies  & notation &  & $f$ verifies  & $g$ verifies  & notation \\
 \hline
 \hline
\multirow{4}{1cm}{Case 2A} &  2A$_1$ &  2A$_1$ & $2A_1-2A_1$ &  \multirow{4}{1cm}{Case 2B}   &  2B$_1$  &  2A$_2$ &$2B_1-2A_2$\\
\cline{2-4}\cline{6-8}
 &  2A$_1$ & 2A$_2$   & $2A_1-2A_2$ &   &2B$_1$  &  2B$_1$ &$2B_1-2B_1$\\
\cline{2-4}\cline{6-8}
 &  2A$_1$ &  2B$_1$ & $2A_1-2B_1$&  & 2B$_1$  &  2B$_2$ &$2B_1-2B_2$\\
\cline{2-4}\cline{6-8}
 &  2A$_1$ &   2B$_2$  &    $2A_1-1B_2$  &    & 2B$_2$  &  2A$_2$ &$2B_2-2A_2$\\
\hline
\end{tabular}
\caption{ Cases to consider for $k=2n$ in the proof of the
stability of $S_\pphapprox$. } \label{table6ch2}
\end{table}

$\bullet$ Case $2A_1-2A_1$: with  equation
(\ref{eq2A1demoth2.7}) we can conclude
\begin{eqnarray*}
|d^2f^1_{2n}-d^2g^1_{2n}| & \leq & \frac{15+3+9\cdot2+5\cdot2}{64}||d^2f-d^2g||_\infty\\
 & \leq &   \frac{23}{32}||d^2f-d^2g||_\infty.\\
\end{eqnarray*}

$\bullet$ Case $2A_1-2A_2$: from  equations
(\ref{eq2A1demoth2.7}) and (\ref{eq2A2demoth2.7}), we obtain
\begin{eqnarray*}
d^2f^1_{2n}-d^2g^1_{2n} &= & \frac{15}{64}d^2f_{n}- \frac{19}{64}d^2g_{n}-\frac{3}{64}d^2f_{n+1}+ \frac{9}{64}\pph(d^2f_{n},d^2f_{n+1})\\
& &  -\frac{3}{64}\pph(d^2g_{n},d^2g_{n+1}) - \frac{5}{64}\pph(d^2f_{n-1},d^2f_{n}) +\frac{7}{64}\pph(d^2g_{n-1},d^2g_{n})) \\
& =& \frac{16}{64}(d^2f_{n}-d^2g_{n})+ \frac{9}{64}(\pph(d^2f_{n},d^2f_{n+1})-\pph(d^2g_{n},d^2g_{n-1}))\\
& &
-\frac{3}{64}(d^2f_{n+1}+d^2g_{n}-2\pph(d^2g_{n},d^2g_{n+1}))\\&&
-
\frac{5}{64}(\pph(d^2f_{n-1},d^2f_{n})-\pph(d^2g_{n-1},d^2g_{n}))\\
& &
-\frac{1}{64}(d^2f_{n}+d^2g_{n-1}-2\pph(d^2g_{n-1},d^2g_{n})).
\end{eqnarray*}
Applying the Lemma \ref{lemme2.1} and the Proposition
\ref{prop2.1}, we obtain
\begin{eqnarray*}
|d^2f^1_{2n}-d^2g^1_{2n}| &\leq  & \frac{16+9\cdot2+3\cdot3+5\cdot2+3}{64}||d^2f-d^2g||_\infty\\
 & \leq & \frac{28}{32}||d^2f-d^2g||_\infty.\\
\end{eqnarray*}

Similarly,

$\bullet$ Case $2A_1-2B_1$: from  equations (\ref{eq2A1demoth2.7})
and (\ref{eq2B1demoth2.7}), we obtain applying the Lemma
\ref{lemme2.1} and the Proposition \ref{prop2.1}
\begin{eqnarray*}
|d^2f^1_{2n}-d^2g^1_{2n}| & \leq & \frac{25}{32}||d^2f-d^2g||_\infty.\\
\end{eqnarray*}

$\bullet$ Case $2A_1-2B_2$: from  equations (\ref{eq2A1demoth2.7})
and (\ref{eq2B2demoth2.7}), we obtain applying the Lemma
\ref{lemme2.1} and the Proposition \ref{prop2.1}
\begin{eqnarray*}
|d^2f^1_{2n}-d^2g^1_{2n}| & \leq & \frac{13}{16}||d^2f-d^2g||_\infty.\\
\end{eqnarray*}

$\bullet$ Case $2B_1-2A_2$: from  equations (\ref{eq2B1demoth2.7})
and (\ref{eq2A2demoth2.7}), we obtain applying the Lemma
\ref{lemme2.1} and the Proposition \ref{prop2.1}
\begin{eqnarray*}
|d^2f^1_{2n}-d^2g^1_{2n}| & \leq & \frac{57}{64}||d^2f-d^2g||_\infty.\\
\end{eqnarray*}

$\bullet$ Case $2B_1-2B_1$: from  equation
(\ref{eq2B1demoth2.7}) we conclude
\begin{eqnarray*}
|d^2f^1_{2n}-d^2g^1_{2n}| & \leq & \frac{16+3+1+9\cdot2+7\cdot2}{64}||d^2f-d^2g||_\infty\\
 & \leq &   \frac{21}{32}||d^2f-d^2g||_\infty.\\
\end{eqnarray*}

$\bullet$ Case $2B_1-2B_2$: from  equations (\ref{eq2B1demoth2.7})
and (\ref{eq2B2demoth2.7}), we obtain applying the Lemma
\ref{lemme2.1} and the Proposition \ref{prop2.1}
\begin{eqnarray*}
|d^2f^1_{2n}-d^2g^1_{2n}| & \leq & \frac{15}{16}||d^2f-d^2g||_\infty.\\
\end{eqnarray*}

$\bullet$ Case $2B_2-2A_2$: from  equations
(\ref{eq2B2demoth2.7}) and (\ref{eq2A2demoth2.7}), we obtain
applying the lemma \ref{lemme2.1} and the proposition
\ref{prop2.1}
\begin{eqnarray*}
|d^2f^1_{2n}-d^2g^1_{2n}| & \leq & \frac{21}{32}||d^2f-d^2g||_\infty.\\
\end{eqnarray*}
\ \\
Finally, the hypotheses of  theorem (\ref{h2th2.3}) are verified
and stability can be established.

$\Box$

\section{Order of approximation}\label{sec5}

In this section we consider the reproduction of polynomials and
the order of a\-ppro\-ximation of $S_{\pphapprox}$.

We recall the following definitions.

\begin{definition}
\label{defreprod}
 A subdivision scheme $S$  is said to reproduce polynomials of  degree $k$
   if for all polynomial $P$ of degree $k$:
\begin{eqnarray*}
if \; \forall n \in \Z, f_n=P(n),\, \textrm{then }  \exists
\tilde{P} \textrm{ a polynomial of degree } k \textrm{ such that }
(Sf)_n=\tilde{P}(2^{-1}n).
\end{eqnarray*}
\end{definition}

\begin{definition}
\label{defordre}
 A subdivision scheme $S$  is said to have an order $k$ of a\-ppro\-xi\-mation if for all function $g \in C^k$ and all $h>0$,
\begin{eqnarray*}
if\; f=g(h.), \; then \,|Sf-g(2^{-1}h.)|\leq C h^k.
\end{eqnarray*}
\end{definition}

We then have the following property.

\begin{proposition}
\label{propo2.7} $S_{\pphapprox}$ reproduces the polynomials of
degree $2$ with translation of $\frac{1}{4}$.
\end{proposition}

{\bf Proof}

We remark that for any $P$,  polynomial of degree $2$, and
$p=(P(n))_{n\in \mathbb{Z}}$, we have
$$
\pph(d^2p_n, d^2p_{n+1})=\frac{p_n+p_{n+1}}{2}.
$$
Therefore, for the initial sequence  $p=(p_n)_{n\in \mathbb{Z}}$,
$S_{\pphapprox}(p)$ coincides with the application to $p$ of  the
linear scheme \cite{DFH05}. In particular, the results of N. Dyn,
M.S. Floater and K. Hormann \cite{DFH05} can be applied and the
property of definition \ref{defreprod} is satisfied with
$\tilde{P}(.)=P(.-1/4)$.

$\Box$

 Concerning the order of approximation the following proposition
 holds.




 \begin{proposition}
\label{prop-ordre} For all function $g\in C^4 ([0,1])$ and $h>0$,
if
 $$f=(g((n-\frac{1}{2})h))_{n\in\mathbb{Z}},$$ then\\
 if $d^2f_nd^2f_{n+1}>0$ for all $n\in \mathbb{Z}$, then
$$||(S_\pphapprox
f)_n-g(2^{-1}h(n-\demi))||_\infty=O(h^4),$$
otherwise
$$||(S_\pphapprox
f)_n-g(2^{-1}h(n-\demi))||_\infty=O(h^3).$$
\end{proposition}

{\bf Proof}

 According to  Proposition \ref{prop2.1}, we have that if
  $d^2f_nd^2f_{n+1}>0$ for all $n\in \mathbb{N}$ then
$$ |\pph(d^2f_n,d^2f_{n+1})-\frac{d^2f_n+d^2f_{n+1}}{2}|=O(h^4).$$
Therefore, is  $S$ stands for the linear scheme  defined in \cite{DFH05},
according to the definition of the $S_{\pphapprox}$,
$$||S_{\pphapprox}f-Sf||_\infty =O(h^4).$$

 Since (see \cite{DFH05}) the scheme $S$ is of order of approximation $4$ we get the result when $d^2f_nd^2f_{n+1}>0$.
In other case, the reproduction of polynomials
leads to
$$||(S_\pphapprox
f)_n-g(2^{-1}h(n-\demi))||_\infty=O(h^3).$$

$\Box$

\begin{remark}

Following \cite{Kuijt98} one can also establish, using the
stability of  $S_\pphapprox$  that $||S^\infty _\pphapprox
f-g||_\infty =O(h^3).$
\end{remark}


\section{Elimination of the Gibbs phenomenon}\label{secG}

\begin{figure}[!ht]

\begin{center}
\begin{tabular}{cc}
\psfig{figure=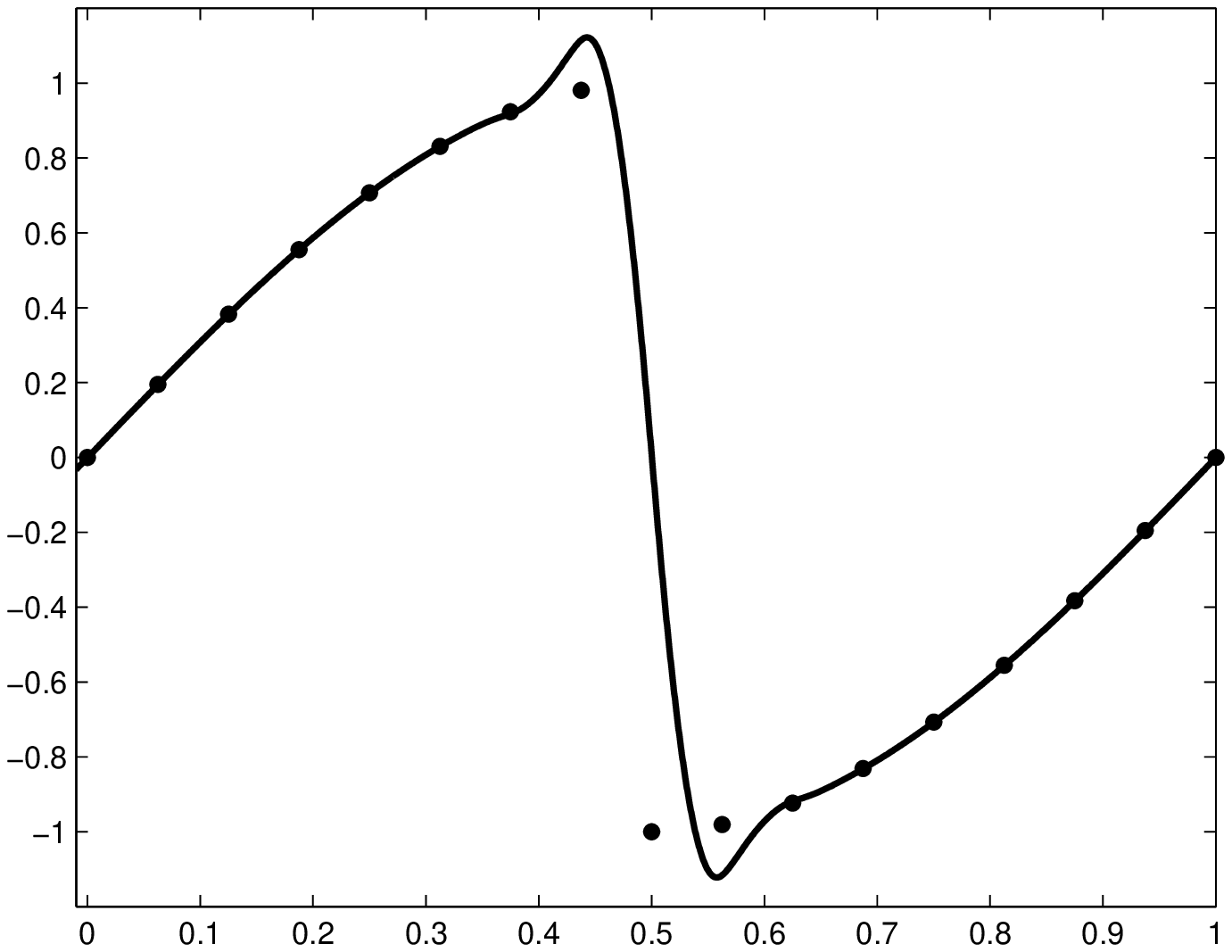,height=6cm}
&\psfig{figure=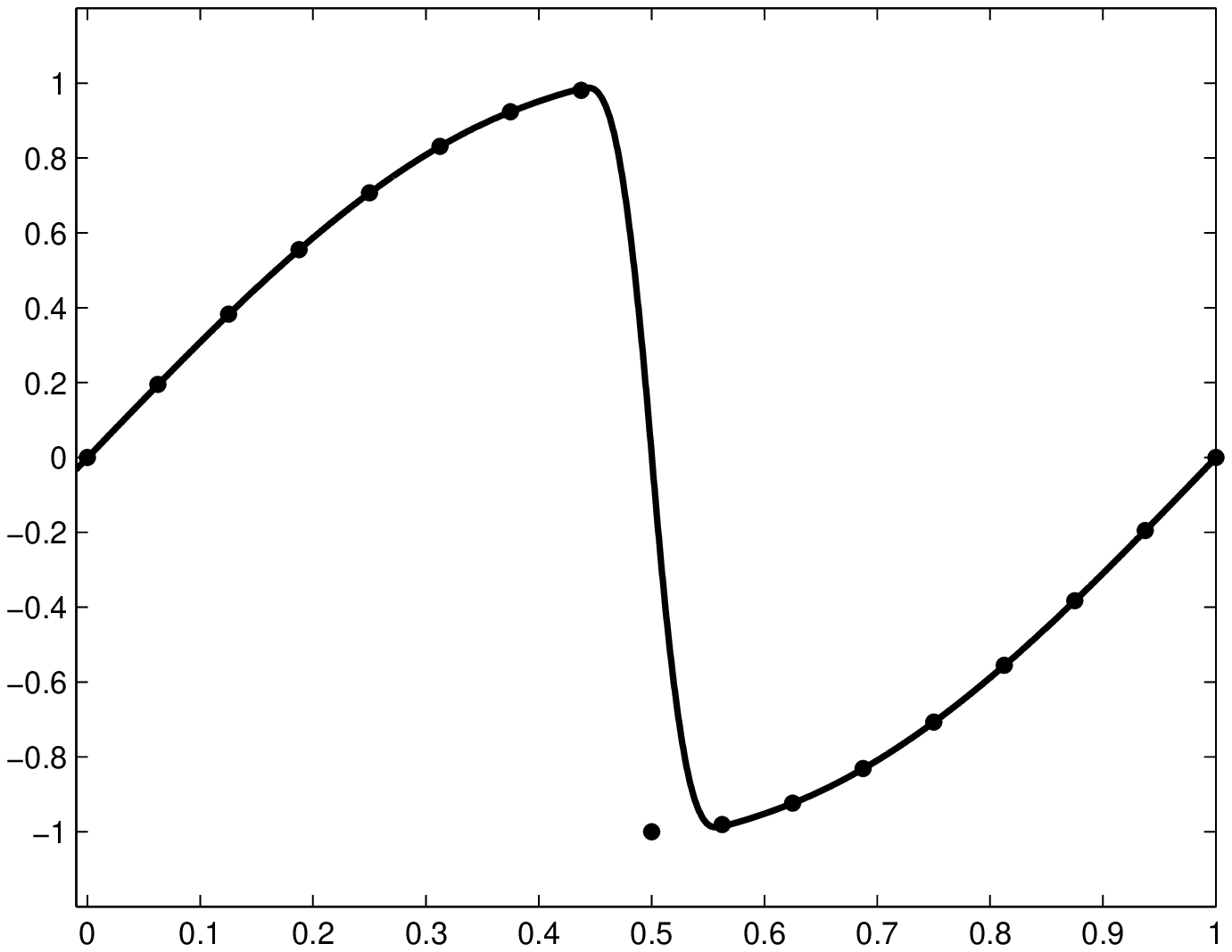,height=6cm}\\
\end{tabular}
\end{center}
 \caption{Comparison of
limit functions for the same initial sequence (sampling of
function (\ref{eqf})).  Left, linear scheme (\ref{schemaDFH}),
 right 
nonlinear scheme $S_\pphapprox$}\label{f1}
\end{figure}

In this section we focus on the behavior of the scheme in presence
of strongly variating data. The reference behavior deals with a
step function as shown on Figure \ref{f1}. As it is visible on
Figure \ref{f1} left, high order linear schemes suffer from an
oscillating behavior named as Gibbs phenomenon.

According to D. Gottlieb and C.W. Shu \cite{GoShu}, given a
punctually discontinuous function $f$ and its sampling $f^h$
defined by $f^h_n=f(nh)$, the Gibbs phenomenon deals with the
convergence of $S^\infty(f^h)$ towards $f$. It can be
characterized by two features (\cite{GoShu} p. 244):
\begin{enumerate}
\item Away from the discontinuity the convergence is rather slow and for any point $x$,
$$|f(x)-S^\infty(f^h)(x)|=O(h).$$
\item There is an overshoot, close to the discontinuity, that does not diminish with reducing  $h$; thus
$$\max|f(x)-S^\infty(f^h)(x)|  \textrm{ does not tend to zero with } h.$$
\end{enumerate}

We are now going to prove that the nonlinear schemes $S_\pphapprox$  does not suffer from the Gibbs phenomenon oscillations,  as it can be guessed from Figure \ref{f1}. We have indeed the following

\begin{proposition}
Given $0\leq \xi\leq h$, for any function $f$  defined by:
\begin{eqnarray*}
\forall x \leq \xi, f(x)=f_-(x)\, \textrm{with} \, f_- \in C^\infty(\left]-\infty, \xi\right], \\
\forall x >\xi, f(x)=f_+(x) \, \textrm{with} \, f_+ \in
C^\infty(\left[\xi, +\infty \right[,
\end{eqnarray*}
and discontinuous in $\xi$,
 we have:

\begin{itemize}
\item if $|x| \geq \frac{9}{2} h,
|f(x-\demi)-S^\infty_\pphapprox(f^h)(x)|=O(h^3),$ \item if $
|x|\leq \frac{9}{2} h , f_-(0)+O(h)\geq
S^\infty_\pphapprox(f^h)(x) \leq f_-(h)+O(h)$.
\end{itemize}

\end{proposition}

{\bf Proof}

Without loss of generality, we focus on $\left[0, +\infty \right[$
and suppose that $f_-(\xi)>f_+(\xi)$.

We first  consider a single application of $S_\pphapprox$. Using
Proposition \ref{prop-ordre} we get:

\begin{itemize}
\item  for $n\geq 2$ and $n_1 \in\left\{2n, 2n+1\right\}$, $|S_\pphapprox(f^h)_{n_1}-f_+(2^{-1}h(n_1-\demi)|=O(h^3)$

\item for $n=1$ since $f$ is discontinuous in $\xi$, $d^2f_n=O(1)$
and $d^2f_{n+1}=O(h^2)$. Then, from  Proposition \ref{prop2.1},
$\pph (d^2f_n, d^2f_{n+1})=O(h^2)$.
 Moreover, according to the definition of  $S_\pphapprox$  as a perturbation of the
 Chaikin scheme $S_C$ we get that $|S_\pphapprox(f^h)_{n_1}-S_C(f^h)_{n_1}|=O(h^2)$.
  Since $S_C$ is a second order scheme  we get that
  $|S_\pphapprox(f^h)_{n_1}-f_+(2^{-1}h(n_1-\demi)|=O(h^2)$ for
  $n_1 \in\left\{2n, 2n+1\right\}$.

\item for $n=0$,  $d^2f_nd^2f_{n+1}\leq 0$ and therefore,
according to the Proposition \ref{prop2.1}, $\pph(d^2f_n,
d^2f_{n+1})=0.$ It is then easy to check, from the definition of
$S_\pphapprox$ that $f_-(0)\leq S_\pphapprox f_{2n}\leq
S_\pphapprox f_{2n+1} \leq f_+(h)$. However, writing
$S_\pphapprox$ as a perturbation of $S_C$ we get that
$|S_\pphapprox(f^h)_{n_1}-S_C(f^h)_{n_1}|=O(d^2f^h)$.

\end{itemize}

Iterating, according to the  stability of $S_\pphapprox$ we get:
\begin{itemize}
\item for $x\geq \frac{9}{2} h$, $|S_\pphapprox^\infty f^h(x) -f_+(x-1/2)|=O(h^3)$.

\item for $ 0\leq x\leq \frac{9}{2} h$, the contraction of the
second order differences (equation (\ref{eq6demoth2.7})) and the
fact that the Chaikin scheme $S_C$ does not produce Gibbs
oscillations  allow to conclude.
\end{itemize}

$\Box$

Before concluding this work, we come back  to  Figure \ref{f1} and
to the comparison between  the limit functions obtained with
$S\pphapprox$ and the limit function obtained  with  linear
subdivision schemes starting from the sampling $f^h$ of the
discontinuous function:

\begin{eqnarray}
\label{eqf}
f(x)&=&\left\{\begin{array}{ll}
sin(\pi x)& \textrm{ for } x \in \left[0,0.5\right]\\
-sin(\pi x)& \textrm{ for } x \in \left]0.5,1\right].
\end{array}
\right.
\end{eqnarray}

It appears from Figure  \ref{f1} that the nonlinear scheme
$S_\pphapprox$ exhibits a much better behavior close to the
discontinuity than the linear scheme of comparable complexity. From Proposition  \ref{prop-ordre} we know moreover that the scheme $S_\pphapprox$ is of higher order than the 
Chaikin scheme.

\section{Conclusions}\label{sec6}

In this paper,  a new nonlinear subdivision scheme has been
defined. It has many   desirable properties. It is  convergent
with  a regularity proved to be at least $C^{1.192-}$ and
numerically estimated at $C^{2.438-}$. By construction, it is
adapted to the presence of isolated discontinuities and  the Gibbs
phenomenon is eliminated. The scheme is  also stable, that due the
nonlinear nature is not a consequence of the convergence.
Moreover, its order of convergence is $3$. Recalling that it is
constructed from a four-point centered stencil, all these
properties make of this scheme a very good candidate for various
applications.


\end{document}